\documentclass[11pt,reqno]{amsart}

\usepackage{amssymb, amsmath, amsthm}
\usepackage{hyperref}
\usepackage[alphabetic,lite]{amsrefs}
\usepackage{verbatim}
\usepackage{amscd}   % for commutative diagrams
\usepackage[all]{xy} % for complicated commutative diagrams
\usepackage{youngtab} % for Young tableaux
\usepackage{young} % for Young tableaux
\usepackage{ytableau}
\usepackage{tikz}
\usepackage{ mathrsfs }
\usepackage{cases}
\usepackage{array}
\usepackage{tabu}
\usepackage{calligra,mathrsfs}

\DeclareMathOperator{\ShHom}{\mathscr{H}\text{\kern -3pt {\calligra\large om}}\,}

\newcommand{\D}{\mathcal{D}}

\newcommand{\DR}{\operatorname{DR}}
\newcommand{\Gr}{\operatorname{Gr}}

\newcommand{\C}{\mathbf{C}}

\renewcommand{\ker}{\operatorname{ker}}

\newtheorem{lemma}{Lemma}[section]

\theoremstyle{definition}

\newtheorem{remark}[lemma]{Remark}
\newtheorem{examplex}[lemma]{Example}
\newenvironment{example}
  {\pushQED{\qed}\examplex}
  {\popQED\endexamplex}

\newtheorem*{main-thm*}{Main Theorem}
\newtheorem*{linear-resolutions*}{Theorem on Linear Resolutions}
\newtheorem*{regularity-powers*}{Theorem on Regularity}
\newtheorem*{injectivity-Ext*}{Theorem on Injectivity of Maps of Ext Modules}
\newtheorem*{Kodaira*}{Kodaira Vanishing for Determinantal Thickenings}

\theoremstyle{definition}

\newtheorem*{definition*}{Definition}

\theoremstyle{remark}
\newtheorem*{remark*}{Remark}

\numberwithin{equation}{section}

\begin{document}

\title{Calculating Hodge and weight filtrations on localizations}

\author{Andr\'as C. L\H{o}rincz}
\address{University of Oklahoma, David and Judi Proctor Department of Mathematics, Norman, OK 73019}
\email{lorincz@ou.edu}

\author{Michael Perlman}
\address{Department of Mathematics, The University of Alabama, Tuscaloosa, AL 35401}
\email{mperlman@ua.edu}

\subjclass[2020]{}

\date{}



\begin{abstract} 
We describe the \textit{Macaulay2} package \texttt{MixedHodgeModules}, which computes Hodge and weight filtrations on localizations $S_f$ of a polynomial ring $S$ at a non-constant reduced polynomial $f$, as well as on twisted localizations $S_f f^{-\alpha}$, where $\alpha\in \mathbf{Q}$ is rational. The package includes functionality for computing Hodge ideals, weighted Hodge ideals, higher multiplier ideals, generation levels of Hodge filtrations, weight levels of elements, $p$-functions, and lengths of weight filtrations. It also provides routines for computing graded de Rham complexes, Du Bois complexes, and intersection Du Bois complexes of hypersurfaces.
\end{abstract}

\maketitle

\section{Introduction}\label{sec:intro}

Let $S=\C[x_1,\cdots,x_n]$ and let $f\in S$ be a reduced non-constant polynomial. The \textit{Macaulay2} \cite{M2} package \texttt{MixedHodgeModules} has functionality for calculating Hodge filtrations, weight filtrations, and related invariants on the localization $S_f$ and, more generally, twisted localizations $S_ff^{-\alpha}$, where $\alpha$ is a rational number. Many of our algorithms are based on \cite{LYang2} and extensions of \cite{blancoArxiv}.

Writing $\D$ for the Weyl algebra of $S$, the $\D$-module $S_ff^{-\alpha}$ underlies a complex mixed Hodge module \cites{saito90,SS}. This implies that it is endowed with two increasing filtrations, indexed by integers: (1) the Hodge filtration $F_{\bullet}$, an infinite filtration by finitely-generated $S$-modules and (2) the weight filtration $W_{\bullet}$, a finite filtration by holonomic $\D$-modules. These two filtrations yield three important classes of ideals in $S$, namely the Hodge ideals \cites{MPHodge, MPQQ}, the weighted Hodge ideals \cites{olanoMultiplier, olanoHodge}, and the higher multiplier ideals \cites{saitoIdeals, SY}.

The Hodge and weight filtrations on $S_ff^{-\alpha}$ are constructed using the graph module $B_f:=\Gamma_+(S)$, obtained via $\D$-module push forward of $S$ along the embedding $\Gamma$ of $\mathbf{A}^n$ whose image is the graph of $f$. The module $B_f$ is endowed with a decreasing filtration $V$, indexed by rational numbers, known as the Kashiwara–Malgrange filtration (see \cite{Vnotes}). We obtain the Hodge-theoretic invariants of $S_ff^{-\alpha}$ discussed above using $V^{\alpha}(B_f)$.

This package calculates the Hodge filtration $F_{\bullet}$ on $V^{\alpha}(B_f)$, as well as kernels of the operators $(s+\alpha)^k$ on $F_p\operatorname{Gr}^{\alpha}_V(B_f)$, using algorithms based on \cite{blancoArxiv}. As applications, this package computes the Hodge ideals, weighted Hodge ideals, and higher multiplier ideals for $\mathbf{Q}$-divisors, as well as related invariants, including the generation level of the Hodge filtration on $S_ff^{-\alpha}$ and Hodge rational homology level \cites{DOR1, PP}.

In another direction, this package has functionality for the $p$-functions of rational functions $(g/f^k)\cdot f^{-\alpha}$ and the related numerical invariants $\nu_{f,g,\alpha}$, based on \cite{LYang2}. As applications, this package can check the weight level of $(g/f^k)\cdot f^{-\alpha}$ in $S_ff^{-\alpha}$, as well as the length of the weight filtration on $S_ff^{-\alpha}$.

Finally, the package uses these Hodge filtration computations to calculate
the Du Bois complexes $\underline{\Omega}^p_{D}$ and the intersection
Du Bois complexes $I\underline{\Omega}^p_{D}$, where $D$ is the divisor associated to $f$.

The routines based on Blanco's algorithm take
$0<\alpha\leq1$. On the other hand, the routines \verb|hodgeCheck|,
\verb|hodgeLevel|, \verb|weightCheck|, \verb|weightLevel|,
\verb|nuAlpha|, and \verb|pFunction| accept every positive
rational $\alpha$. We specify the parameter ranges below. These restrictions concern the implementation, rather than
the definition of twisted localizations.

\subsection*{Organization} Section~\ref{sec:hodge}  recalls the graph module $B_f$ and describes calculations of Hodge filtrations on twisted localizations, Hodge ideals, higher multiplier ideals, generation levels, and Hodge rational homology levels. Section~\ref{sec:weight} discusses the weight filtration on $S_f f^{-\alpha}$, weighted Hodge ideals, weight levels of elements, $p$-functions, and the length of the weight filtration. Section~\ref{sec:complexes} describes the routines for graded de Rham complexes, Du Bois complexes, and intersection Du Bois complexes.

\section{The Hodge filtration and related invariants}\label{sec:hodge}

Let $\D=S\langle \partial_1,\cdots,\partial_n \rangle$ be the Weyl algebra of $S$. Our convention is that all $\D$-modules are left $\D$-modules. We write $F_{\bullet}(\D)$ for the increasing filtration on $\D$ by order of differential operator.

Let $f\in S$ be a reduced non-constant function, and let $\alpha\in \mathbf{Q}$. We consider the $\D$ module $S_ff^{-\alpha}$, which has elements $hf^{-\alpha}$ where $h\in S_f$ and $f^{-\alpha}$ is a formal symbol. The $S$-action is given by
\begin{equation}
    g\cdot hf^{-\alpha}=ghf^{-\alpha},\quad \textnormal{for $g\in S$,}
\end{equation}
so that $S_ff^{-\alpha}$ is isomorphic as an $S$-module to the usual localization $S_f$. The action of a derivation $\Theta\in \operatorname{Der}_{\C}(S)$ is given by
\begin{equation}
 \Theta\cdot hf^{-\alpha}=\left(\Theta(h) - \alpha h \frac{\Theta(f)}{f}\right)f^{-\alpha} .  
\end{equation}
Since $S$ underlies a pure Hodge module, $S_ff^{-\alpha}$ is functorially endowed with the structure of a complex mixed Hodge module \cite{SS}. For $\alpha\in\mathbf Q$, write
$\alpha=\beta+k$, where $0<\beta\leq1$ and $k\in\mathbf Z$.
The identification
\begin{equation}\label{eq:twisting}
S_ff^{-\alpha}\xrightarrow{\sim}S_ff^{-\beta},
\qquad hf^{-\alpha}\longmapsto (h/f^k)f^{-\beta},
\quad h\in S_f,
\end{equation}
preserves the Hodge and weight filtrations.
Unless otherwise indicated, the computational routines
discussed in this section take $0<\alpha\leq1$.

Let $\Gamma:\mathbf{A}^n\to \mathbf{A}^n\times \mathbf{A}^1_t$ be the map satisfying $\Gamma(p)=(p,f(p))$. We write $B_f=\Gamma_+(S)$ for the $\D$-module push-forward of $S$ along $\Gamma$. This is a module over $\D\langle t, \partial_t\rangle$, the Weyl algebra on $\mathbf{A}^n\times \mathbf{A}^1_t$. 

If we write $\delta$ for the class of $1/(f-t)$ in $B_f$, then we have
\begin{equation}\label{eqn:graph}
B_f=\bigoplus_{j\geq 0}S\,\partial_t^j\delta.  
\end{equation}
The actions of $\Theta\in \operatorname{Der}_{\C}(S)$ and $t$ are given by
\begin{equation}
\Theta\cdot h\partial_t^j\delta=\Theta(h)\partial_t^j\delta-h\Theta(f)\partial_t^{j+1}\delta,\quad \quad \textnormal{and}\quad\quad t\cdot h\partial_t^j\delta = hf\partial_t^j\delta-jh\partial_t^{j-1}\delta,    
\end{equation}
where $h$ is in $S$. The actions of $S$ and $\partial_t$ are straightforward from (\ref{eqn:graph}). 

\subsection{Hodge filtration on the $V$-filtration}

The Hodge filtration $F_{\bullet}$ on $B_f$ is an increasing $\mathbf{Z}$-indexed filtration of $S[t]$-modules given by $F_{p}(B_f)=0$ for $p<0$ and
\begin{equation}
F_p(B_f)= \bigoplus_{j= 0}^p S\partial_t^j\delta,\quad \textnormal{for $p\geq 0$}.   
\end{equation}
The Kashiwara--Malgrange filtration on $B_f$ is a decreasing filtration $V^{\bullet}$ by $S[t]$-modules, indexed by rational numbers. It is exhaustive, discrete, and left-continuous. See \cite{Vnotes} for a general introduction. For $\alpha\in \mathbf{Q}$ we write $V^{\alpha}(B_f)$ for the corresponding piece of the $V$-filtration, and we write
\begin{equation}
V^{>\alpha}(B_f)= \bigcup_{\beta>\alpha} V^{\beta}(B_f),\quad\quad\textnormal{and}\quad\quad \operatorname{Gr}_V^{\alpha}(B_f)=V^{\alpha}(B_f)/V^{>\alpha}(B_f). 
\end{equation}
The Hodge filtration on $B_f$ induces filtrations $F_{\bullet}$ on $V^{\alpha}(B_f)$ and $\operatorname{Gr}_V^{\alpha}(B_f)$ which we also call Hodge filtrations.

\begin{remark}
Our convention for the Hodge filtration on $B_f$ differs from the usual Hodge filtration by a shift of one, i.e. our \(F_p(B_f)\) is \(F^{\mathrm{std}}_{p+1}(B_f)\). On the other hand, Hodge filtrations on twisted localizations $S_ff^{-\alpha}$ retain their usual indexing.    
\end{remark}

For $0<\beta\leq1$, $k\in\mathbf Z_{\geq0}$, and $p\geq0$,
multiplication by $t^k$ induces an isomorphism
\[
t^k:F_pV^\beta(B_f)\xrightarrow{\sim}F_pV^{\beta+k}(B_f).
\]
See \cite[Proposition 3.19(i)]{Vnotes}.
Thus, positive indices reduce to $(0,1]$.

We first discuss functionality for $F_p(V^{\alpha}(B_f))$ for $\alpha\in (0,1]$. By discreteness and continuity of the $V$-filtration, there exist rational numbers $1=\alpha_1>\alpha_2>\cdots>\alpha_k=0$ such that $F_p(V^{\alpha}(B_f))$ is constant for $\alpha$ in intervals of the form $(\alpha_{i+1},\alpha_i]$. In other words, 
\[
F_p(V^{\alpha}(B_f))=F_p(V^{\beta}(B_f))\quad \iff \quad \alpha,\beta\in (\alpha_{i+1},\alpha_i] \textnormal{ for some $1\leq i<k$}.
\]
The set $\{\alpha_1,\cdots,\alpha_{k-1}\}$ can be described in terms of roots of the $(p+1)$-st generalized $b$-function $b^{(p+1)}_f(s)$, which is the monic polynomial of minimal degree satisfying the functional equation
\begin{equation}
b^{(p+1)}_f(s)\delta = P(s)t^{p+1}\delta,\quad \textnormal{for some $P(s)\in \D[s]$,}    
\end{equation}(see \cite[Section 5]{BL} and \cite[Algorithm 1]{blancoArxiv}). Note that these roots are in fact determined by the roots of the usual Bernstein-Sato polynomial $b_f(s):= b_f^{(1)}(s)$, but not their multiplicities. More precisely, 
\[
\alpha\in \{\alpha_1,\cdots,\alpha_{k-1}\}\quad \implies \quad b_f^{(p+1)}(-\alpha-p)=0 \textnormal{ and } \alpha\in (0,1].
\]
The function \verb|hodgeOnV| takes as input $f$ and $p$, and outputs a hash table whose keys are $\alpha\in(0,1]$ such that $-\alpha$ is a root of $b_f^{(p+1)}(s-p)$, and whose values are lists representing $S$-generators for $F_p(V^{\alpha}(B_f))$.

\begin{example}
We consider the $\mathsf{D}_4$ singularity $f=x^2+y^3+yz^2\in \mathbf{Q}[x,y,z]$ and set $p=1$.

\vspace{.1cm}

\begin{verbatim}
i1 : needsPackage "MixedHodgeModules"
i2 : S = QQ[x,y,z];
i3 : f = x^2+y^3+y*z^2;
i4 : hodgeOnV(f,1)
                          3                 2      2
o4 = HashTable{1 => {1, -z dt, -y*z*dt, - 3y dt - z dt, -x*dt}}
               1
               - => {1, -z*dt, -y*dt, -x*dt}
               2
               1
               - => {1, -dt}
               6
               5          2               2
               - => {1, -z dt, -y*z*dt, -y dt, -x*dt}
               6               
\end{verbatim}

\vspace{.1cm}

\noindent We see $\{\alpha_1,\alpha_2,\alpha_3,\alpha_4\}=\{1,\frac{5}{6},\frac{1}{2},\frac{1}{6}\}$ are the keys of the hash table. The value associated to the key $\alpha_i$ is an $S$-generating set for $F_1(V^{\alpha_i}(B_f))$. We compare to the roots of $b^{(2)}_f(s)$ using the function \verb|generalB| in the \verb|BernsteinSato| package of \textit{Macaulay2}:

\vspace{.1cm}

\begin{verbatim}
i5 : factorBFunction generalB({f},1_S, Exponent => 2)    
                        3      5      7      11      13      17
o5 = (s + 1)(s + 2)(s + -)(s + -)(s + -)(s + --)(s + --)(s + --)
                        2      2      6       6       6       6
\end{verbatim}

\vspace{.1cm}

\noindent We see that $\{1,\frac{5}{6},\frac{1}{2},\frac{1}{6}\}=\{2-1,\frac{11}{6}-1,\frac{3}{2}-1,\frac{7}{6}-1\}$.

Alternatively, the user may also input some $\alpha\in (0,1]$ to get an $S$-generating set of $F_p(V^{\alpha}(B_f))$:

\vspace{.1cm}

\begin{verbatim}
i6 : hodgeOnV(f,3/4,1)
           2               2
o6 = {1, -z dt, -y*z*dt, -y dt, -x*dt}
\end{verbatim}

\vspace{.1cm}

\noindent Comparing outputs, we see that $F_1(V^{\frac{3}{4}}(B_f))=F_1(V^{\frac{5}{6}}(B_f))$, which is expected, as $\frac{3}{4}\in (\frac{1}{2},\frac{5}{6}]$.
\end{example}

Our algorithm is based on \cite[Algorithm 1]{blancoArxiv}. We remark that the published version \cite{blanco} of this algorithm has a minor issue, which has subsequently been corrected in the arXiv preprint \cite{blancoArxiv}.

\subsection{Hodge filtrations on twisted localizations} Next, we consider the Hodge filtration on $S_ff^{-\alpha}$. For $j\geq 0$ let $Q_j(x)\in \mathbf{C}[x]$ be the polynomial $Q_j(x)=\prod_{i=0}^{j-1}(x+i)$, with the convention that $Q_0(x)=1$. Given a rational number $\alpha\in (0,1]$, define
\begin{equation}
\tau_{\alpha}:V^{\alpha}(B_f)\to S_ff^{-\alpha},\quad\quad \textnormal{via}\quad\quad \tau_{\alpha}\left(\sum_{j=0}^p h_j\partial_t^j\delta \right) =\sum_{j=0}^p Q_j(\alpha)\frac{h_j}{f^{j+\alpha}}.  
\end{equation}
This is a surjection for all $\alpha\in (0,1]$. The Hodge filtration $F_{\bullet}$ on $S_ff^{-\alpha}$ is given by 
\begin{equation}
F_p(S_ff^{-\alpha})= \tau_{\alpha}(F_p(V^{\alpha}(B_f))),\quad \textnormal{for $p\geq 0$},       
\end{equation}
see \cite[Corollary~1.2]{DY25}. We say that the Hodge filtration on $S_ff^{-\alpha}$ is generated in level $\ell$ if
\[
F_p(\D)\cdot F_{\ell}(S_ff^{-\alpha})= F_{\ell+p}(S_ff^{-\alpha})\quad \textnormal{for all $p\geq 0$},
\]
where $F_p(\D)$ denotes the differential operators of order at most $p$. The generation level of the Hodge filtration on $S_ff^{-\alpha}$ is the minimal $\ell$ such that the Hodge filtration is generated in level $\ell$. If $n\geq 2$ and $\alpha=1$, then the generation level is $\leq n-2$, and the generation level is $\leq n-1$ otherwise (see \cite[Theorem E]{MPlocal} for a stronger statement involving the minimal exponent).

The functions \verb|doesGenerateNext| and \verb|generationLevel| test generation of the Hodge filtration. The function \verb|doesGenerateNext(f,alpha,p)| checks whether $F_1(\D)$ applied to $F_p(S_ff^{-\alpha})$ generates $F_{p+1}(S_ff^{-\alpha})$.

\begin{example}
For the $2\times 2$ determinant, the Hodge filtration on $S_f$ has generation level one.

\vspace{.1cm}

\begin{verbatim}
i1 : S = QQ[x,y,z,w];
i2 : f = x*w-y*z;

i3 : doesGenerateNext(f,1,0)
o3 = false

i4 : doesGenerateNext(f,1,1)
o4 = true

i5 : generationLevel(f)
o5 = 1
\end{verbatim}

\vspace{.1cm}

\noindent For the cusp, the generation level is zero when $\alpha=1$.

\vspace{.1cm}

\begin{verbatim}
i6 : S = QQ[x,y];
i7 : f = x^2+y^3;

i8 : doesGenerateNext(f,1,0)
o8 = true

i9 : generationLevel(f)
o9 = 0
\end{verbatim}
\end{example}

\subsection{Hodge ideals} We write $D$ for the reduced divisor associated to $f\in S$, and we let $\alpha\in (0,1]$. For $p\geq 0$ the Hodge ideal $I_p(\alpha D)\subseteq S$ is the ideal defined by the following equality \cites{MPHodge, MPQQ}:
\begin{equation}
I_p(\alpha D)f^{-p-\alpha}=F_p(S_ff^{-\alpha}),\quad \textnormal{for $p\geq 0$}.    
\end{equation}
The function \verb|hodgeIdeal| computes the Hodge ideal $I_p(\alpha D)$. The user inputs \verb|hodgeIdeal(f,alpha,p)|, where \verb|f| is the defining equation of $D$. The general generation level bounds discussed in the previous section are used for efficiency. As a first sanity check, we calculate for a normal crossing divisor.

\begin{example}
Let $f=xyz\in \mathbf{Q}[x,y,z,w]$. We calculate the first few Hodge ideals for $\alpha=1$:

\vspace{.1cm}

\begin{verbatim}
i2 : S = QQ[x,y,z,w];
i3 : f = x*y*z;
i4 : alpha = 1;

i5 : hodgeIdeal(f,alpha,0)
o5 = ideal 1

i6 : hodgeIdeal(f,alpha,1)
o6 = ideal (y*z, x*z, x*y)

i7 : hodgeIdeal(f,alpha,2)
             2 2       2   2 2     2    2      2 2
o7 = ideal (y z , x*y*z , x z , x*y z, x y*z, x y )

i8 : hodgeIdeal(f,alpha,3)
             3 3     2 3   2   3   3 3     3 2   2 2 2   3   2   2 3    3 2    3 3
o8 = ideal (y z , x*y z , x y*z , x z , x*y z , x y z , x y*z , x y z, x y z, x y )
\end{verbatim}

\vspace{.1cm}

\noindent This is consistent with \cite[Proposition 8.2]{MPHodge} in the case $n=4$ and $r=3$.
\end{example}

In the next example, we observe $\alpha$ dependence of $I_p(\alpha D)$.

\begin{example}\label{e:sym}
Let $f\in \mathbf{Q}[x_1,\cdots,x_6]$ be the determinant of the $3\times 3$ symmetric matrix. When $\alpha=1$:  

\vspace{.1cm}

\begin{verbatim}
i2 : S = QQ[x_1..x_6];
i3 : f = determinant genericSymmetricMatrix(S,x_1,3);

i4 : alpha = 1;

i5 : hodgeIdeal(f,alpha,0)
o5 = ideal 1

i6 : hodgeIdeal(f,alpha,1)
             2                                    2                       2
o6 = ideal (x  - x x , x x  - x x , x x  - x x , x  - x x , x x  - x x , x  - x x )
             5    4 6   3 5    2 6   3 4    2 5   3    1 6   2 3    1 5   2    1 4          
             \end{verbatim}

\vspace{.1cm}

\noindent On the other hand, when $\alpha=1/2$ we obtain a different result:

\vspace{.1cm}

\begin{verbatim}
i7 : alpha = 1/2;

i8 : hodgeIdeal(f,alpha,0)
o8 = ideal 1

i9: hodgeIdeal(f,alpha,1)
o9 = ideal 1
\end{verbatim}

\vspace{.1cm}

\noindent Compare to \cite[Corollary 5.4]{LYangSemi}, which describes Hodge ideals for all $\mathbf{Q}$-divisors associated to symmetric determinants.
\end{example}

The package also includes a direct implementation of Zhang's formula for Hodge ideals of $\mathbf{Q}$-divisors with weighted homogeneous isolated singularities \cite{zhang}. The function \verb|hodgeIdealWeightedHomogIsolated| takes as input a weighted homogeneous polynomial $f$, a rational number $\alpha\in (0,1]$, a Hodge level $p$, and a list of weights $w_1,\ldots,w_n$ normalized so that every monomial of $f$ has weight one.

\begin{example}
We consider the $\mathsf{E}_8$ singularity $f=x^2+y^3+z^5\in \mathbf{Q}[x,y,z]$, with weight list $w=\{1/2,1/3,1/5\}$.

\vspace{.1cm}

\begin{verbatim}
i2 : S = QQ[x,y,z];
i3 : f = x^2+y^3+z^5;
i4 : w = {1/2,1/3,1/5};
i5 : alpha = 1/2;

i6 : hodgeIdealWeightedHomogIsolated(f,alpha,1,w)
                     2   3
o6 = ideal (x, y*z, y , z )

i7 : hodgeIdealWeightedHomogIsolated(f,alpha,2,w)  
                      2   2    3     3   2 2    3      2    4   5   3    2   2 3   4    5
o7 = ideal (x*y*z, x*y , x y, x , x*z , x z , 3y z - 2x z, y , z + y - 2x , y z , y z, y )
\end{verbatim}
\end{example}

The functions \verb|hodgeCheck(f,g,alpha,p)| and
\verb|hodgeLevel(f,g,alpha)| respectively test membership of
$gf^{-\alpha}$ in $F_p(S_ff^{-\alpha})$ and determine its
Hodge level. Both accept every positive rational $\alpha$,
using the identification (\ref{eq:twisting}).

\subsection{Higher multiplier ideals} 

For $p\geq 0$ the higher multiplier ideal $\tilde{I}_p(\alpha D)\subseteq S$ is defined via the equality \cites{saitoIdeals, SY}:
\begin{equation}
\tilde{I}_p(\alpha D)\otimes \partial_t^p=\operatorname{Gr}^F_p(V^{\alpha}(B_f)),\quad \textnormal{for $p\geq 0$}.    
\end{equation}
We always have (see \cites{MPQQ, SY}):
\[
I_p(\alpha D)+(f)=\tilde{I}_p(\alpha D)+(f).
\]

The function \verb|higherMultiplierIdeal| computes the higher multiplier ideal $\widetilde{I}_p(\alpha D)$.

\begin{example}
For the cusp $f=x^2+y^3$, we recover \cite[Equation (1.8)]{DY25}.

\vspace{.1cm}

\begin{verbatim}
i1 : S = QQ[x,y];
i2 : f = x^2+y^3;
i3 : alpha = 11/12;

i4 : higherMultiplierIdeal(f,alpha,2)
             2    3     3   5
o4 = ideal (x y, x , x*y , y )

i5 : hodgeIdeal(f,alpha,2)
             3    4      2      3   2 2
o5 = ideal (x , 6y  - 17x y, x*y , x y )
\end{verbatim}
\end{example}

\begin{example}
The next calculation is consistent with \cite[Proposition D]{zhang} and \cite[Conjecture E]{zhang}:

\vspace{.1cm}

\begin{verbatim}
i6 : S = QQ[x,y];
i7 : f = x^2+y^3;
i8 : alpha = 1;
i9 : w = {1/2,1/3};

i10 : higherMultiplierIdeal(f,alpha,0)
o10 = ideal (y, x)
i11 : hodgeIdealWeightedHomogIsolated(f,alpha,0,w)
o11 = ideal (y, x)


i12 : higherMultiplierIdeal(f,alpha,1)
                   2   3
o12 = ideal (x*y, x , y )
i13 : hodgeIdealWeightedHomogIsolated(f,alpha,1,w)
                   2   3
o13 = ideal (x*y, x , y )

i14 : higherMultiplierIdeal(f,alpha,2)
              2    3     3   5
o14 = ideal (x y, x , x*y , y )
i15 : hodgeIdealWeightedHomogIsolated(f,alpha,2,w)
              3   4     2      3   2 2
o15 = ideal (x , y  - 3x y, x*y , x y )
\end{verbatim}
\end{example}

\subsection{The Hodge rational homology level}
We continue to write $D=V(f)$. The local cohomology Hodge module $H^1_D(S)=S_f/S$ satisfies
\[
W_{n+1}(H^1_D(S))=\operatorname{IC}_{D}^H(-1),
\]
where $\operatorname{IC}^H_{D}(-1)$ is the intersection cohomology Hodge module of $D$ with a Tate twist of $-1$, and $W_{\bullet}$ denotes the weight filtration on $H^1_D(S)$. The Hodge rational homology level \cites{DOR1, PP}, written $\operatorname{HRH}(D)$, is an invariant that compares the Hodge filtration on $H^1_D(S)$ to that on $\operatorname{IC}_{D}^H(-1)$. More precisely,
\[
\operatorname{HRH}(D)\geq k\quad \iff \quad F_p(H^1_D(S))=F_p(\operatorname{IC}_{D}^H(-1))\textnormal{ for all $p\leq k$,}
\]
under the natural inclusion $\operatorname{IC}_D^H(-1)\subseteq H^1_D(S)$. The property $\operatorname{HRH}(D)\geq k$ is referred to as the condition $D_k$ in \cite{PP}, where it is shown that it implies certain symmetries of the Hodge--Du Bois diamond of $D$. It is known that, if $D$ has rational singularities, then $\operatorname{HRH}(D)\geq 0$ \cites{olanoMultiplier}. We say that $D$ is a rational homology manifold if $H^1_D(S)=\operatorname{IC}_{D}^H(-1)$ (equivalently, $\operatorname{HRH}(D)\geq k$ for all $k\geq 0$). By convention, $\operatorname{HRH}(D)=-1$ if $F_0(H^1_D(S))\neq F_0(\operatorname{IC}_{D}^H(-1))$.

We implement this invariant in \textit{Macaulay2} using \cite[Theorem A]{DOR2}. The function \verb|HRHCheck(f,k)| tests whether $\operatorname{HRH}(f)\geq k$, while \verb|HRHLevel(f)| returns the maximal $k$ for which $\operatorname{HRH}(D)\geq k$, or the string \verb|rational homology manifold| when appropriate.

\begin{example}
The $3\times 3$ generic determinant has HRH level zero.

\vspace{.1cm}

\begin{verbatim}
i2 : S = QQ[x_(1,1)..x_(3,3)];
i3 : f = determinant genericMatrix(S,x_(1,1),3,3)

i4 : HRHCheck(f,0)
o4 = true

i5 : HRHCheck(f,1)
o5 = false

i6 : HRHLevel(f)
o6 = 0
\end{verbatim}

\vspace{.1cm}

\noindent The property $\operatorname{HRH}(D)\geq 0$ is consistent with the fact that the $3\times 3$ determinant has rational singularities. 
\end{example}

All generic determinantal varieties have HRH level zero, see \cites{PRaicu, PHodge}.

\begin{example}
In contrast, the $\mathsf{A}_1$ singularity is a rational homology manifold.

\vspace{.1cm}

\begin{verbatim}
i6 : S = QQ[x,y,z];
i7 : f = y^2-x*z;

i8 : HRHCheck(f,0)
o8 = true

i9 : HRHCheck(f,1)
o9 = true

i10 : HRHLevel(f)
o10 = rational homology manifold
\end{verbatim}
\vspace{.1cm}

All finite quotient singularities are rational homology manifolds.
\end{example}

\section{The weight filtration and related invariants}\label{sec:weight}

In this section, we discuss weighted Hodge ideals, adjoint ideals, $p$-functions, and weight length of $S_f f^{-\alpha}$. 

\subsection{Monodromy weight filtration}

Let $s=-\partial_tt$ and $\alpha\in \mathbf{Q}_{>0}$. Multiplication by $N=(s+\alpha)$ is  nilpotent on $\operatorname{Gr}_V^{\alpha}(B_f)$, and we have the monodromy weight filtration $W(N)_{\bullet}$ on $\operatorname{Gr}_V^{\alpha}(B_f)$ centered at $0$:
\begin{equation}\label{eqn:monodromyWeight}
W(N)_m    =\sum_{i+j=m} \ker(N^{i+1})\cap \operatorname{im}(N^{-j}),\quad \textnormal{for $m\in \mathbf{Z}$},    
\end{equation}
where $\ker(N^k)=0$ for $k\leq 0$ and $\operatorname{im}(N^k)=\operatorname{Gr}^{\alpha}_V(B_f)$ for $k\leq 0$. This is an increasing filtration by $\D$-modules. 

Due to algorithmic limitations, we are unable to implement the full filtration $W(N)_{\bullet}$. However, we can still get partial information, sufficient for calculating the weighted Hodge ideals (see Section \ref{weightOnLoc}).

\subsection{The weight filtration on twisted localizations}\label{weightOnLoc}

The ring $S$ is a pure Hodge module of weight $n$, and the (complex) mixed Hodge module $S_ff^{-\alpha}$ has weights $\geq n$. Continuing to write $W(N)_{\bullet}$ for the monodromy weight filtration on $\operatorname{Gr}^{\alpha}_V(B_f)$ centered at zero, for $m\geq 0$ we consider the lift of $W(N)_{m-1}$ to $V^{\alpha}(B_f)$:
\begin{equation}
W_{m-1}(V^{\alpha}(B_f)):=\{v\in V^{\alpha}(B_f)\mid \overline{v}\in W(N)_{m-1}\}.   
\end{equation}
The weight filtration on $S_ff^{-\alpha}$ is given by (e.g. \cite[Section 4.2]{LYang2})
\begin{equation}
W_{m+n}(S_ff^{-\alpha})=\tau_{\alpha}(W_{m-1}(V^{\alpha}(B_f))),\quad \textnormal{for $m\geq 0$}. 
\end{equation}
The weight length of $S_ff^{-\alpha}$ is defined to be the largest $m\geq 0$ such that $\operatorname{Gr}^W_{m+n}(S_ff^{-\alpha})$ is nonzero.

We can calculate $W_{m+n}(S_ff^{-\alpha})$ using a submodule of $W_{m-1}(V^{\alpha}(B_f))$. For $m\geq 0$, let $K_m^{\alpha}$ be the lift of $\operatorname{ker}(N^m)$ to $V^{\alpha}(B_f)$, i.e.
\begin{equation}
K_m^{\alpha}=\{v\in V^{\alpha}(B_f)\mid N^m\overline{v}=0\;\textnormal{in}\; \operatorname{Gr}^{\alpha}_V(B_f)\},    
\end{equation}
so $K_m^{\alpha}\subseteq W_{m-1}(V^{\alpha}(B_f))$ using (\ref{eqn:monodromyWeight}). By \cite[Section 4.2]{LYang2} we have  (see also \cite[Theorem A]{olanoHodge}):  
\begin{equation}
    W_{m+n}(S_ff^{-\alpha})=\tau_{\alpha}(K_m^{\alpha}).
\end{equation}
 
\iffalse

\begin{lemma}
Using notation as above, we have
\[
W_{m+n}(S_ff^{-\alpha})=\tau_{\alpha}(K_m^{\alpha}).
\]
\end{lemma}

\begin{proof}
By (\ref{eqn:monodromyWeight}) we have $\operatorname{ker}(N^m)\subseteq W(N)_{m-1}$, so $K_m^{\alpha}\subseteq W_{m-1}(V^{\alpha}(B_f))$, and thus $\tau_{\alpha}(K_m^{\alpha})\subseteq W_{m+n}(S_ff^{-\alpha})$.

For the reverse containment, we note that $W(N)_{m-1}\subseteq \ker(N^m)+\operatorname{im}(N)$.
The preimage of $\operatorname{im}(N)\subseteq \operatorname{Gr}_V^{\alpha}(B_f)$ in $V^{\alpha}(B_f)$ is
\[
(s+\alpha)V^{\alpha}(B_f)+V^{>\alpha}(B_f),
\]
and $V^{>\alpha}(B_f)\subseteq K^{\alpha}_m$. Thus, it suffices to show that $\tau_{\alpha}((s+\alpha)v)=0$
for all $v=\sum_{j=0}^q h_j\partial_t^j\delta\in V^{\alpha}(B_f)$. We compute
\[
(s+\alpha)v=\sum_{j=0}^q(\alpha+j)h_j\partial_t^j\delta-h_jf\partial_t^{j+1}\delta.
\]
Applying $\tau_{\alpha}$ gives
\[
\tau_{\alpha}((s+\alpha)v)
=\sum_{j=0}^q\left((\alpha+j)Q_j(\alpha)-Q_{j+1}(\alpha)\right)\frac{h_j}{f^{j+\alpha}}=0,
\]
since $Q_{j+1}(\alpha)=(\alpha+j)Q_j(\alpha)$. Therefore $\tau_{\alpha}((s+\alpha)v)=0$, as claimed.
\end{proof}

\fi

The modules $K_m^{\alpha}$ can be calculated using the function \verb|weightHodgeOnV|, which uses a refinement of Blanco's algorithm \cite{blancoArxiv}. This function facilitates the calculation of the weighted Hodge ideals, discussed in the next sub-section. We remark that the modules $K_m^{\alpha}$ have a description via multiplicities of roots of Bernstein-Sato polynomials, similar to Sabbah's description of $V$-filtrations, see \cite[Section 4.2]{LYang2}. The reason that we cannot calculate the full monodromy weight filtration $W(N)_{\bullet}$ is that Blanco's algorithm is not amenable to calculating images $\operatorname{im}(N^k)$. 

\subsection{Weighted Hodge ideals}\label{sec:weightedIdeals}

For $m,p\geq 0$, the weighted Hodge ideal $I^{W_m}_p(\alpha D)\subseteq S$ is defined via the equality \cites{olanoMultiplier, olanoHodge}:
\begin{equation}\label{def:weightedHodge}
I_p^{W_m}(\alpha D)f^{-p-\alpha}=F_pW_{m+n}(S_ff^{-\alpha}),\quad \textnormal{for $m,p\geq 0$}.
\end{equation}
We remark that \cite{olanoHodge} only defines weighted Hodge ideals for reduced divisors ($\alpha=1$). The definition (\ref{def:weightedHodge}) is a straightforward generalization of the one in loc. cit. For a fixed $p$ and $\alpha$, we get an increasing chain
\[
I_p^{W_0}(\alpha D)\subseteq I_p^{W_1}(\alpha D)\subseteq I_p^{W_2}(\alpha D)\subseteq \cdots \subseteq I_p(\alpha D)\subseteq S,
\]
where $I_p^{W_m}(\alpha D)=I_p(\alpha D)$ for sufficiently large $m$ (at most the weight length of $S_ff^{-\alpha}$).

The function \verb|weightedHodgeIdeal| computes the weighted Hodge ideals $I_p^{W_m}(\alpha D)$, for $\alpha\in (0,1]$.

\begin{example}
Let $f\in \mathbf{Q}[x_1,\cdots,x_6]$ be the determinant of the $3\times 3$ symmetric matrix. By \cite{LYangSemi} we have $\operatorname{Gr}^W_p(S_f)=0$ if $p\neq 6,7,8$, and
\[
\operatorname{Gr}^W_6(S_f)=S,\quad \operatorname{Gr}^W_7(S_f)=\operatorname{IC}_D,\quad \operatorname{Gr}^W_8(S_f)=\operatorname{IC}_{\{0\}}.
\]
This is reflected in the following weighted Hodge ideal calculation:

\vspace{.1cm}

\begin{verbatim}
i2 : S = QQ[x_1..x_6];
i3 : f = determinant genericSymmetricMatrix(S,x_1,3);
i4 : alpha = 1;
i5 : weightedHodgeIdeal(f,alpha,1,1)
             2                                    2                       2
o5 = ideal (x  - x x , x x  - x x , x x  - x x , x  - x x , x x  - x x , x  - x x )
             5    4 6   3 5    2 6   3 4    2 5   3    1 6   2 3    1 5   2    1 4 
\end{verbatim}

\vspace{.1cm}

\noindent Comparing to Example \ref{e:sym} above, we see that $I_1(D)=I_1^{W_1}(D)$. On the other hand, with the following code

\vspace{.1cm}

\begin{verbatim}
i6 : weightedHodgeIdeal(f,alpha,2,1)
i7 : weightedHodgeIdeal(f,alpha,2,2)
\end{verbatim}

\vspace{.1cm}

\noindent we obtain that $I_2^{W_1}(D)$ has $21$ minimal generators, each of degree four, whereas $I_2^{W_2}(D)$ has $16$ minimal generators, $15$ of degree four, and one of degree three (corresponding to $f^{-2}$ in $S_f$). Therefore, $I_2^{W_1}(D)\subsetneq I_2^{W_2}(D)$, which indicates that the Hodge filtration on $\operatorname{IC}_{\{0\}}$ induced by that on $S_f$ starts in level two. See also \cite[Corollary 5.4]{LYangSemi}.
\end{example}

\begin{example}
On the other hand, when $f=x_1^2+\cdots+x_5^2\in \mathbf{Q}[x_1,\cdots,x_5]$, the divisor $D$ is a rational homology manifold:

\vspace{.1cm}

\begin{verbatim}
i2 : S = QQ[x_1..x_5];
i3 : f = x_1^2+x_2^2+x_3^2+x_4^2+x_5^2;
i4 : HRHLevel(f)
o4 = rational homology manifold
\end{verbatim}

\vspace{.1cm}

\noindent In particular, $I_p(D)=I_p^{W_1}(D)$ for all $p\geq 0$:

\vspace{.1cm}

\begin{verbatim}
i5 : alpha = 1;
i6 : hodgeIdeal(f,alpha,1)
o6 = ideal 1
i7 : weightedHodgeIdeal(f,alpha,1,1)
o7 = ideal 1
i8 : hodgeIdeal(f,alpha,2)
o8 = ideal (x , x , x , x , x )
             5   4   3   2   1
i9 : weightedHodgeIdeal(f,alpha,2,1)
o9 = ideal (x , x , x , x , x )
             5   4   3   2   1
\end{verbatim}

\vspace{.1cm}

We see that $I_1(D)=I_1^{W_1}(D)$ and $I_2(D)=I_2^{W_1}(D)$.
\end{example}

The adjoint ideal is the weighted Hodge ideal for $p=0$ and $\alpha=1$, namely $\operatorname{adj}(D)=I^{W_1}_0(D)$ (see \cite[Theorem A]{olanoMultiplier} and \cite[Section 9.3E]{Laz}). The divisor $D$ has rational singularities if and only if $\operatorname{adj}(D)=S$. The function \verb|adjointIdeal| computes this ideal directly.

\subsection{Weight level}

For $\alpha \in \mathbf{Q}$ and $g \in S$, denote by $\omega_{f,g,\alpha}$ the (shifted) weight level of the element $g f^{-\alpha}$, that is, the smallest integer $\omega \geq 0$ for which we have $g f^{-\alpha} \in W_{n +\omega} (S_f f^{-\alpha})$. We implement the algorithm in \cite[Algorithm 5.9]{LYang2} for the determination of these invariants through the function $ \verb|weightLevel(f,g,alpha)|$.

\begin{example}\label{KashBar}
We consider the cubic cone
$f=x^3+y^3+z^3\in\mathbf{Q}[x,y,z]$ and take
$g=\partial_x f=3x^2$. We compare the weight levels of
$g/f$ and $1/f$.

\vspace{.1cm}

\begin{verbatim}
i2 : S = QQ[x,y,z];
i3 : f = x^3+y^3+z^3;
i4 : g = diff(x,f);

i5 : weightLevel(f,g,1)
o5 = 1

i6 : weightLevel(f,1_S,1)
o6 = 2
\end{verbatim}
The first output is consistent with the result in \cite{barletkashiwara} that implies that $\omega_{f,\frac{\partial f}{\partial x},1}=1$ for any irreducible polynomial $f$ with $\frac{\partial f}{\partial x}\neq 0 .$
\end{example}

This function relies on another function, \verb|weightCheck|, which takes $f$, $g$, $\alpha$, $\ell$ and checks if $gf^{-\alpha}$ belongs to $W_{n+\ell}(S_ff^{-\alpha})$. We emphasize that, unlike \verb|hodgeCheck| and \verb|hodgeLevel|, the functions \verb|weightCheck| and \verb|weightLevel| do not require knowledge of an entire piece of the filtration. Rather, the weight functionality relies on \cite{LYang2}. 

\subsection{Power $b$-functions and $p$-functions}
For an integer $k\geq 1$, the $k$-th power (generalized) Bernstein-Sato polynomial $b^{(k)}_{f,g}(s)$ of a function $g$ with respect
     to $f$ is given by the minimal monic polynomial in $\mathbf{Q}[s]$ satisfying an equation
     \begin{equation}\label{eqn: function eqn for higher b function} 
P\cdot g f^{k}\delta=b_{f,g} ^{(k)}(s)g \delta.
\end{equation}
     with $P(s) \in \D[s]$. This agrees, after a rescaling,
     with the usual generalized Bernstein-Sato polynomial of $g$ with respect to $f^k$. When $k=1$, we set $b_{f,g}(s):= b_{f,g}^{(k)}(s)$, and have $b_f(s)= b_{f,1}(s)$ (the usual Bernstein-Sato polynomial of $f$).

While the roots of $b_{f,g}(s)$ determine those of $b_{f,g}^{(k)}(s)$, it does not determine the multiplicities of the roots. As in \cite{LYang2}, for $\alpha \in \mathbf{Q}$ we define
\[\nu_{f,g, \alpha}:= \lim_{k\to \infty}\operatorname{mult}_{s=-\alpha}b^{(k)}_{f,g}(s).\]
The limit stabilizes for an explicit bound on $k$ \cite[Theorem 1.3]{LYang2}. By \emph{loc. cit.}, the numbers $\nu_{f,g, \alpha}$ provide an upper bound for the weight level, that is $\omega_{f,g,\alpha} \leq \nu_{f,g,\alpha}$, so that
\[ g f^{-\alpha} \in W_{n+\nu_{f,g, \alpha}}(S_f f^{-\alpha}).\]

The function \verb|nuAlpha(f,g,alpha)| computes these numbers $\nu_{f,g, \alpha}$. 

\begin{example}
Continuing Example \ref{KashBar}, we have

\begin{verbatim}
i7 : nuAlpha(f,g,1)
o7 : 1

i8 : nuAlpha(f,1_S,1)
o8 : 2
\end{verbatim}
Note that the first output is consistent with the general fact from \cite{LYang2} that $\nu_{f,\frac{\partial f}{\partial x},1}= \omega_{f,\frac{\partial f}{\partial x},1}=1$ for any irreducible polynomial $f$ with $\frac{\partial f}{\partial x}\neq 0 $ . The second output is consistent with the Torelli problem, as well as the general fact that $\nu_{f, 1, \alpha} = \omega_{f, 1, \alpha}$ for isolated quasi-homogeneous singularities (see \cite{LYang2}).

The default strategy computes using the annihilating ideal in $\D[s]$ of $\delta$. Another strategy, \verb|PowerBFunction|, obtains $\nu$ by computing the Bernstein-Sato polynomial of a suitable power of $f$. A third strategy, implemented only when $g=1$, is called \verb|Malgrange| and circumvents the computation of $\operatorname{Ann}_{\D[s]}(\delta)$ using the Malgrange ideal (see \cite[(5.20)]{SST} for the definition of this ideal). The user toggles the strategy using \verb|NuMethod =>|. For instance,
\begin{verbatim}
i9 : nuAlpha(f,1_S,1, NuMethod => Malgrange)
o9 : 2    
\end{verbatim}
\end{example}

Following \cite{LYangSemi, LYang2}, the $p$-function $p_{f,g,\alpha}(s)$ is the monic polynomial of minimal degree in $\mathbf{Q}[s]$ such that $p_{f,g,\alpha}(s) g \delta \in V^{\alpha}(B_f)$, where $V^\bullet$ again denotes the Kashiwara-Malgrange $V$-filtration of $B_f$ with respect to $f$. By \cite{LYang2}, the $p$-function can be expressed as a finite product
\[ p_{f,g,\alpha}(s)=\prod_{\beta<\alpha}(s+\beta)^{\nu_{f,g, \beta}}.\]
By \emph{loc. cit.} the degree of $p_{f,g,\alpha}(s)$ gives an upper bound on the Hodge level of $g f^{-\alpha}$, that is
     \[ g f^{-\alpha} \in F_{\deg p_{f,g,\alpha}(s)}(S_f f^{-\alpha}).\]
The function \verb|pFunction(f,g,alpha)| computes the polynomials $p_{f,g, \alpha}(s)$. The calculation is based on computing $\nu$, for which one can use the strategies \verb|NuMethod =>| as described above.

\begin{example}
Let $f\in\mathbf{Q}[x_1,\cdots,x_6]$ be the determinant of the
$3\times3$ symmetric matrix, and let
$g=x_1x_4-x_2^2=\partial f/\partial x_6$ be its upper-left
$2\times2$ minor. 
\vspace{.1cm}

\begin{verbatim}
i2 : S = QQ[x_1..x_6];
i3 : f = determinant genericSymmetricMatrix(S,x_1,3);
i4 : g = x_1*x_4-x_2^2;

i5 : pFunction(f,g,2)
o5 = s + 1

i6 : factorBFunction pFunction(f,1_S,2)
                 3
o6 = (s + 1)(s + -)
                 2
\end{verbatim}
\vspace{.1cm}

\noindent The polynomial $g$ can be realized as a highest-weight vector of weight
$(2,2,0)$ under the natural action of the $3\times 3$ general linear group. By \cite[Proposition 5.4]{LYangSemi},
\[
b_{f,g}(s)=(s+1)(s+\tfrac52)(s+3),\qquad
b_f(s)=(s+1)(s+\tfrac32)(s+2).
\]
The outputs agree with \cite[Theorem 1.1]{LYangSemi}. Moreover,
\cite[Theorem 1.4]{LYangSemi} shows that their degrees give the
Hodge levels exactly:
\[
\frac{g}{f^2}\in F_1(S_f)\setminus F_0(S_f),\qquad
\frac{1}{f^2}\in F_2(S_f)\setminus F_1(S_f).
\]
For the twist $\alpha=\tfrac12$, we obtain:

\vspace{.1cm}

\begin{verbatim}
i7 : pFunction(f,g,3/2)
o7 = s + 1
\end{verbatim}

\vspace{.1cm}

\noindent Again using \cite[Theorem 1.4]{LYangSemi}, we obtain
\[
\frac{g}{f}f^{-1/2}\in
F_1(S_ff^{-1/2})\setminus F_0(S_ff^{-1/2}).
\]
\end{example}

\subsection{Length of the weight filtration} For $\alpha \in \mathbf{Q}$, we call the weight length the smallest $n_{f,\alpha} \in \mathbf{Z}_{\geq 0}$ such that $W_{n+n_{f,\alpha}}(S_f f^{-\alpha}) = S_f f^{-\alpha}$. This agrees with the nilpotency index of the logarithmic monodromy operator $N=(s+\alpha)$ on the nearby cycles $\operatorname{Gr}^{\alpha}_V (B_f)$.

Note that for any $g\in S$ we have $\nu_{f,g,\alpha} \leq n_{f, \alpha}$ (see \cite{LYang2}). An interesting fact is that $\omega_{f,g, \alpha} = n_{f,\alpha}$ if and only if $\nu_{f,g, \alpha} = n_{f,\alpha}$. In \cite[Corollary 1.9]{LYang2}, a result is given for determining $n_{f,\alpha}$ via the Bernstein-Sato polynomial of a certain power $k_0$ of $f$.

The function \verb|weightLength(f,alpha)| calculates the number $n_{f,\alpha}$. The default strategy performs a calculation of $\nu$. Alternatively, as for the same power $k_0$ the weight level is equal to the corresponding $\nu$, one can opt for the optional strategy \verb|ByWeightLevel| using \verb|LengthStrategy =>|. The default strategy is \verb|ByNuAlpha|, in which case the user can select sub-strategies \verb|NuMethod =>| as described in the previous subsection.

\begin{example}
Let $S=\mathbf{Q}[x,y,z]$ and let $f=xyz(x+y+z)$. We first compute the Bernstein--Sato polynomial.

\vspace{.1cm}

\begin{verbatim}
i2 : S = QQ[x,y,z];
i3 : f = x*y*z*(x+y+z);

i4 : factorBFunction globalBFunction f
            3     3      3      5
o4 = (s + 1) (s + -)(s + -)(s + -)
                  2      4      4

i5 : weightLength(f,1)
o5 = 3

i6 : weightLength(f,1,LengthStrategy => ByWeightLevel)
o6 = 3
\end{verbatim}

\vspace{.1cm}

\noindent Thus $n_{f,1}=3$, or equivalently,
\[
W_5(S_f)\subsetneq W_6(S_f)=S_f.
\]

Following \cite[Corollary 1.9]{LYang2}, we test the twists corresponding
to $-\tfrac14,-\tfrac12,-\tfrac34$. We also test $\alpha=\tfrac13$, for which
$-\alpha+\mathbf{Z}$ contains no root of $b_f(s)$.

\vspace{.1cm}

\begin{verbatim}
i7 : apply({1/4,1/2,3/4}, alpha -> weightLength(f,alpha))
o7 = {1, 1, 1}

i8 : weightLength(f,1/3)
o8 = 0
\end{verbatim}

\vspace{.1cm}

\noindent Hence
\[
n_{f,1/4}=n_{f,1/2}=n_{f,3/4}=1,\qquad n_{f,1/3}=0.
\]
For $\alpha\in\{\tfrac14,\tfrac12,\tfrac34\}$, this means
\[
W_3(S_ff^{-\alpha})\subsetneq
W_4(S_ff^{-\alpha})=S_ff^{-\alpha}.
\]
In contrast, $S_ff^{-1/3}$ is pure of weight three.
\end{example}

\section{Graded de Rham and Du Bois complexes}\label{sec:complexes}
In this section, we consider complexes of finitely generated $S$-modules arising from the Hodge filtrations on $H^1_D(S)$ and $\operatorname{IC}_D^H$. Let $(M,F_{\bullet})$ be the filtered $\D$-module underlying a mixed Hodge module $\mathcal{M}$. The de Rham complex $\DR_S(\mathcal{M})$ is the $\C$-linear complex
\begin{equation}
0\to M\to \Omega^1_S\otimes_S M \to \cdots \to \omega_S\otimes_S M\to 0,    
\end{equation}
lying in cohomological degrees $-n,\cdots ,0$. This complex has a filtration $F_{\bullet}$ induced by the one on $M$, with associated graded complexes $\operatorname{Gr}^F_p\DR_S(\mathcal{M})$:
\begin{equation}
0\to \Gr^F_p(M) \to \Omega^1_S\otimes_S \Gr^F_{p+1}(M)\to \cdots \to \omega_S\otimes_S \Gr^F_{p+n}(M)\to 0,\quad p\in \mathbf{Z}.   
\end{equation}
These are complexes of finitely generated $S$-modules, in contrast to the complex $\DR_S(\mathcal{M})$.

\subsection{Graded de Rham complexes}

The function \verb|gradedDeRhamComplexH1| computes a free complex quasi-isomorphic to $\operatorname{Gr}^F_p(\operatorname{DR}(H^1_D(S)))$, for $p\geq -n$. As complexes in \textit{Macaulay2}
     are homologically graded, the $-i$-th homology of this complex is the $i$-th
     cohomology of $\operatorname{Gr}^F_p(\operatorname{DR}(H^1_D(S)))$.

\begin{example}
We calculate several graded pieces for $f=x^2+y^3+z^2\in \mathbf{Q}[x,y,z]$.

\vspace{.1cm}

\begin{verbatim}
i2 : S = QQ[x,y,z];
i3 : f = x^2+y^3+z^2;
i4 : C1 = gradedDeRhamComplexH1(f,-3)
      1      1
o4 = S  <-- S          
     0      1

i5 : prune HH_0(C1) 
o5 = cokernel | y3+x2+z2 |

i6 : C2 = gradedDeRhamComplexH1(f,-2)
      3      10      7
o6 = S  <-- S   <-- S
     0      1       2
     
i7 : prune HH_1(C2)
o7 = cokernel | -y  0  x  -z  |
              | x   z  y2 0   |
              | -3z 3x 0  3y2 |
              | 0   y  -z -x  |

i8 : C3 = gradedDeRhamComplexH1(f,-1)
      6      27      37      16
o8 = S  <-- S   <-- S   <-- S
     0      1       2       3
\end{verbatim}

\vspace{.1cm}

In particular, $\mathcal{H}^{-1}(\operatorname{Gr}^F_{-2}(\operatorname{DR}(H^1_D(S))))\neq 0$.
\end{example}

The function \verb|gradedDeRhamCohomologyH1| computes a specified cohomology module without constructing the entire complex.

\subsection{Du Bois complexes}

The function \verb|duBoisComplex| computes a free complex quasi-isomorphic to the Du Bois complex $\underline{\Omega}^p_D$ of $D=V(f)$, which may be identified as (see \cite[Proposition 13.1]{MPLC})
\begin{equation}
\underline{\Omega}_D^p=\mathbf{R}\operatorname{Hom}_S(\Gr^F_{p-n}\DR_S(H^1_D(S)), \omega_S)[p+1].    
\end{equation}
We calculate $\mathbf{R}\operatorname{Hom}_S$ using the \verb|Complexes| package of \textit{Macaulay2}. Again, as complexes in \textit{Macaulay2}
     are homologically graded, the $-i$-th homology of this complex is the $i$-th
     cohomology of $\underline{\Omega}^p_D$. Thus, the homology of the output complex vanishes in positive degree. 

\begin{example}
If $D$ is nonsingular, the Du Bois complexes recover the K\"{a}hler differentials of $D$.

\vspace{.1cm}

\begin{verbatim}
i2 : S = QQ[x,y,z];
i3 : f = y^2-x;

i4 : DB0 = duBoisComplex(f,0)
      1      1
o4 = S  <-- S
     0      1

i5 : prune HH_0(DB0)
o5 = cokernel | y2-x |

i6 : DB1 = duBoisComplex(f,1)
      5      8      3
o6 = S  <-- S  <-- S
     0      1      2
     
i7 : prune HH_0(DB1)
o7 = cokernel | y2-x 0    |
              | 0    y2-x |

i8 : DB2 = duBoisComplex(f,2)
      10      25      21      6
o8 = S   <-- S   <-- S   <-- S
     0       1       2       3
     
i9 : prune HH_0(DB2)
o9 = cokernel | y2-x |
\end{verbatim}

\smallskip

\noindent In particular, $\mathcal{H}^0(\underline{\Omega}^p_D)$ is a free $S/(f)$-module of rank $\binom{2}{p}$. One can show that all other cohomology vanishes.
\end{example}

We write $\Omega^p_D$ for the $p$-th module of K\"{a}hler differentials of $D$. For $p\geq 0$, there is a map $\Omega^p_D\to \underline{\Omega}^p_D$, and if such maps are quasi-isomorphisms for all $p\leq m$, we say that $D$ has $m$-Du Bois singularities \cites{MOPW,SaitoDB}. The classical notion of Du Bois (or log canonical) singularities coincides with $0$-Du Bois. We say $D$ has pre $m$-Du Bois singularities if $\mathcal{H}^i(\underline{\Omega}^p_D)=0$ for all $i>0$ and $p\leq m$ \cite{SVV}, which is a necessary condition for $m$-Du Bois singularities. This property can be checked using the function \verb|isPreDuBois| in the package.

 \begin{example}
The surface $D=V(x^2+y^3+z^7)$ is not pre $0$-Du Bois, since $\mathcal{H}^1(\underline{\Omega}^0_D)$ is nonzero.

 \vspace{.1cm}

 \begin{verbatim}
i2 : S = QQ[x,y,z];
i3 : f = x^2+y^3+z^7;

i4 : DB0 = duBoisComplex(f,0)

       1      4      3
 o4 = S  <-- S  <-- S
      -1     0      1

i5 : prune HH_(-1)(DB0)
o5 = cokernel | z y x |

i6 : isPreDuBois(f,0)
o6 = false   
\end{verbatim}

\smallskip

\noindent We see that $\mathcal{H}^1(\underline{\Omega}^0_D)\cong S/(x,y,z)$, which is supported on the singular locus of $D$. 
\end{example}

\subsection{Intersection Du Bois complexes} Let $\operatorname{IC}_D^H$ be the intersection cohomology Hodge module of $D$. In the hypersurface case, one has
\begin{equation}
\operatorname{IC}^H_D(-1)=W_{n+1}(H^1_D(S)).    
\end{equation}
The intersection Du Bois complex of $D$ is 
\begin{equation}
I \underline{\Omega}_D^p=\Gr^F_{-p}\DR_S(\operatorname{IC}_D^H)[p+1-n].
\end{equation}
Alternatively, since $\operatorname{IC}^H_D$ is self-dual up to a Tate twist, this complex is also represented as
\[
I \underline{\Omega}_D^p=\mathbf{R}\operatorname{Hom}_S(\Gr^F_{p+1-n}\DR_S(\operatorname{IC}^H_D), \omega_S)[p+1].
\]
The function \verb|intersectionDuBoisComplex| computes a free complex quasi-isomorphic to the intersection Du Bois complex $I\underline{\Omega}^p_D$. The homology of the output complex vanishes in positive degree.

\begin{example}
Since the $\mathsf{A}_1$ singularity is a rational homology manifold, the Du Bois and intersection Du Bois complexes agree. We see this for $p=1$:

\vspace{.1cm}

\begin{verbatim}
i2 : S = QQ[x,y,z];
i3 : f = y^2-x*z;
i4 : DB1 = duBoisComplex(f,1)
      7      10      3
o4 = S  <-- S   <-- S
     0      1       2

i5 : IDB1 = intersectionDuBoisComplex(f,1)
      7      10      3
o5 = S  <-- S   <-- S
     0      1       2

i6 : prune HH_0(DB1)
o6 = cokernel | z y 0 -z |
              | y x 0 -y |
              | y x z 0  |
              | 0 0 y x  |

i7 : prune HH_0(IDB1)
o7 = cokernel | z y 0 -z |
              | y x 0 -y |
              | y x z 0  |
              | 0 0 y x  |
\end{verbatim}
\end{example}

\begin{example}
For the $2\times 2$ determinant, we consider $I\underline{\Omega}^p_D$. Since $D$ has rational singularities and admits a small resolution, $\mathcal{H}^0(\underline{\Omega}^p_D)\cong\mathcal{H}^0(I\underline{\Omega}^p_D)\cong\Omega^{[p]}_D$ \cite{KS}, where the latter denotes the reflexive hull of the K\"{a}hler differentials $\Omega^{p}_D$.

%\vspace{.1cm}

\begin{verbatim}
i2 : S = QQ[x,y,z,w];
i3 : f = x*w-y*z;
i4 : IDB0 = intersectionDuBoisComplex(f,0)
      1      1
o4 = S  <-- S
     0      1

i5 : prune HH_0(IDB0)
o5 = cokernel | yz-xw |

i6 : IDB1 = intersectionDuBoisComplex(f,1)
      1      12      15      4
o6 = S  <-- S   <-- S   <-- S
     -1     0       1       2

i7 : prune HH_0(IDB1)
o7 = cokernel | z  yz-xw 0     z2  xz  |
              | w  0     0     zw  yz  |
              | y  0     0     xw  xy  |
              | -x 0     yz-xw -xz -x2 |

i8 : prune HH_(-1)(IDB1)
o8 = cokernel | w z y x |

i9 : IDB2 = intersectionDuBoisComplex(f,2)
      4      49      90      55      10
o9 = S  <-- S   <-- S   <-- S   <-- S
     -1     0       1       2       3

i10 : prune HH_0(IDB2)
o10 = cokernel | x  0 y  0 w y 0  0  |
               | -z 0 -w 0 0 0 0  0  |
               | -y w 0  y 0 0 0  0  |
               | 0  z y  x 0 0 0  0  |
               | 0  0 0  0 w y -w -z |
               | 0  0 0  0 z x 0  0  |
               | 0  0 0  0 0 0 y  x  |

i11 : prune HH_(-1)(IDB2)
o11 = 0
\end{verbatim}

%\vspace{.1cm}
We remark that, using $\mathcal{H}^{n-p-1}(I\underline{\Omega}^p_D)=\mathcal{H}^0(\operatorname{Gr}^F_{-p-1}\operatorname{DR}_S(\operatorname{IC}_D^H(-1)))$, one may calculate the generation level of the Hodge filtration on $\operatorname{IC}_D^H(-1)$, see \cite[Lemma 10.1]{MPLC}. More precisely, the Hodge filtration on $\operatorname{IC}_D^H(-1)$ is generated in level $k$ if and only if $\mathcal{H}^{n-p-1}(I\underline{\Omega}^p_D)=0$ for $p<n-k-1$.
\end{example}

\section*{Acknowledgments}
We thank Guillem Blanco, Bradley Dirks, Timothy Duff, Mahrud Sayrafi, and Ruijie Yang for helpful conversations. This work began as part of the \textit{Macaulay2 Workshop and Mini-School} held at University of Minnesota - Twin Cities, funded by NSF Award DMS-2302476. Perlman acknowledges the support of NSF Award DMS-2601624.

\end{document}